\documentclass{article}
\usepackage{latexsym}
\usepackage{graphics}
\usepackage{amsfonts}
\newtheorem{thm}{Theorem}

\newtheorem{prop}{Proposition}
\newtheorem{lem}{Lemma}
\newtheorem{defn}{Definition}
\newtheorem{cor}{Corollary}
\newtheorem{quest}{Question}

\begin{document}

\title{Tree Diagrams for String Links II:  Determining Chord Diagrams}
\author{Blake Mellor\\
			Mathematics Department\\
			Loyola Marymount University\\
			Los Angeles, CA  90045-2659\\
   {\it  bmellor@lmu.edu}}
\date{}
\maketitle

\begin{abstract}
In previous work \cite{me3}, we defined the intersection graph of a chord diagram
associated with a string link (as in the theory of finite type invariants).  In this paper,
we look at the case when this graph is a tree, and we show that in many cases
these trees determine the chord diagram (modulo the usual 1-term and 4-term relations).
\end{abstract}
\tableofcontents

\section{Introduction} \label{S:intro}

The theory of finite type invariants allows us, via the Kontsevich
integral \cite{ko}, to identify knot and link invariants with their associated
{\it weight systems}, functionals on {\it chord diagrams} which obey
certain relations.  Since these diagrams are purely combinatorial
objects, we can take a combinatorial approach to studying these
weight systems, providing a new viewpoint on the associated knot
invariants.  One approach to these chord diagrams for knots, due to Chmutov,
Duzhin and Lando \cite{cdl}, is to study their {\it intersection
graphs}.  While these graphs do not contain all of the information of
the chord diagrams, they distinguish chord diagrams in many cases
\cite{cdl, me}, and have an interesting algebraic structure of their
own \cite{la}.  Recently, the author has extended the notion of an intersection graph to
string links \cite{me3}.

In this paper we will not discuss the background of finite type invariants, instead looking at the
relationship between the intersection graph and the chord diagram from a purely
combinatorial viewpoint.  For a discussion of how these diagrams arise in the theory of
finite type invariants, see \cite{bn}.  In section~\ref{S:prelim} we will review the
definitions of chord diagrams and intersection graphs for string links.  In section~\ref{S:2comp} we
will look at chord diagrams on two components.  In this case, we look at a special class of trees ({\it
trimmed} trees) which arise as intersection graphs, and show that they determine the
associated chord diagram modulo some standard relations, the 1-term and 4-term relations.  In
section~\ref{S:>2comp} we address the (easier) case of diagrams with more than two components. 
Finally, in section~\ref{S:questions} we pose some questions for further research.

{\it Acknowledgement:}  The author thanks Loyola Marymount University
for supporting this work via a Summer Research Grant in 2004.

\section{Preliminaries} \label{S:prelim}

\subsection{Chord Diagrams} \label{SS:chord}

We begin by defining what we mean by a {\it chord diagram}.  Since we are only considering
chord diagrams which arise from string links, we will simply refer to these as
chord diagrams, but the reader should be aware that these differ somewhat from the more
usual chord diagrams for knots.

\begin{defn}
A {\bf chord diagram of degree n with k components} is a disjoint union of k oriented
line segments (called the {\bf components} of the diagram), together with $n$ chords (unoriented line
segments with endpoints on the components), such that all of the $2n$ endpoints of the chords are
distinct.  The diagram is determined by the orders of the endpoints on each component.
\end{defn}

\begin{figure}
    (1-term relation) $$\includegraphics{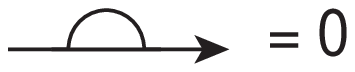}$$
    (4-term relation) $$\includegraphics{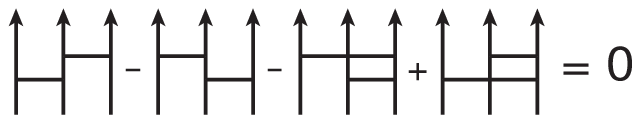}$$
    (antisymmetry relation) $$\includegraphics{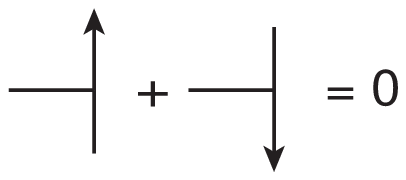}$$
    \caption{The 1-term, 4-term and antisymmetry relations.} \label{F:4-term}
    \end{figure}

We can naturally organize the diagrams with $k$ components into a graded vector space
with real coefficients, graded by their degree.  We denote the vector space of chord
diagrams of degree $n$ on $k$ components by $B_n^k$.  We impose three relations on
$B_n^k$ (motivated by knot theory), called the {\it 1-term}, {\it 4-term} and {\it antisymmetry}
relations, shown in Figure~\ref{F:4-term} (no other chords have endpoints on the arcs shown; in the
4-term relation, all other chords of the four diagrams are the same); we will still call the resulting
vector space $B_n^k$.  The three arcs in the 4-term relation may belong to the same component or to
different components.

It is sometimes useful to combine the vector spaces $B_n^k$ into a graded module
$B^k = \bigoplus_{n\geq 1}B_n^k$ via direct sum.  We can give the module
$B^k$ a bialgebra (or Hopf algebra) structure for any $k$ by defining an appropriate
product and co-product:
\begin{itemize}
    \item  We define the (noncommutative) product $D_1 \cdot D_2$ of two chord diagrams
$D_1$ and $D_2$ as the result of placing $D_2$ on top of $D_1$ (joining the components
so the orientations agree), as shown below:
$$\includegraphics{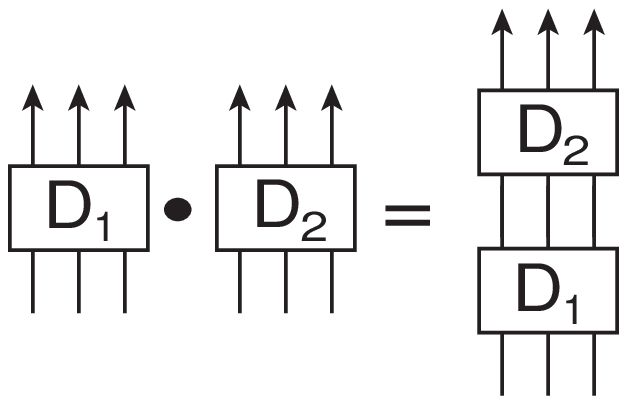}$$
    \item  We define the co-product $\Delta(D)$ of a chord diagram $D$ as
follows:
$$\Delta(D) = {\sum_J D_J' \otimes D_J''}$$
where $J$ is a subset of the set of chords of $D$, $D_J'$ is $D$
with all the chords in $J$ removed, and $D_J''$ is $D$ with all
the chords {\it not} in $J$ removed.  For example:
$$\includegraphics{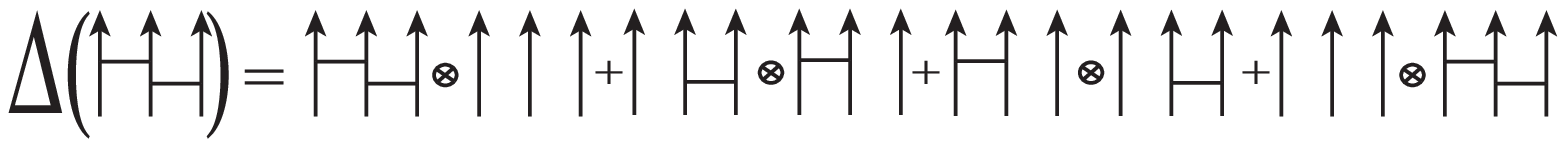}$$
\end{itemize}
It is easy to check the compatibility condition $\Delta(D_1\cdot
D_2) = \Delta (D_1)\cdot\Delta(D_2)$.

We can also define an action of $B^1$ on $B^k$ (in fact, $k$ different actions).  Bar-Natan \cite{bn}
showed that $B^1$ is isomorphic to the space of chord diagrams on knots, modulo the 4-term relation. 
He also showed that there is a well-defined (again, modulo the 4-term relation) commutative product on
this space, where $D_1 \# D_2$ is the connected sum of the diagrams $D_1$ and $D_2$.  Extending these
results, we define the product $A \#_i D$ (where $A \in B^1$ and $D \in B^k$) to be the result of
taking the connected sum of $A$ (viewed as a chord diagram on a circle) and the $i$th component of
$D$.  This action is well-defined modulo the 4-term relation; i.e. it does not matter where on the
$i$th component of $D$ we "glue in" $A$.  For example:
$$\includegraphics{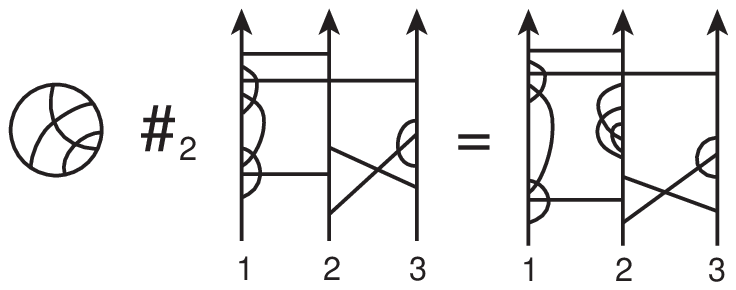}$$

To end this section, we will define the notion of a {\it share} - a collection of chords in a diagram
which can sometimes be treated as a single chord (see, for example, Lemma \ref{L:gen4term} in section
\ref{SS:IGC}).  Shares will be an extremely useful tool in our proofs.

\begin{defn} \label{D:share}
A {\bf share} of a chord diagram D is a subset S of the set of chords in D and two arcs A
and B on the boundary components of D (A and B may be on the same or different components)
such that every chord in S has both endpoints in $A \cup B$ and {\bf no} other chord in D
has an endpoint in $A \cup B$.
\end{defn}

\subsection{Intersection Graphs} \label{SS:IG}

The essential value of the intersection graph for knots (in which the chord diagram
consists of chords in a bounding circle) is that it can detect when the order of two
endpoints for different chords along the bounding circle is switched, since this changes the
pair of chords from (visually) intersecting to non-intersecting or vice-versa.  For chord
diagrams for string links, the existence of a "bottom" and "top" for each component allows us
to give a linear (rather than cyclic) ordering to the endpoints of the chords on each
component, and so the notion of one endpoint being "below" another is well-defined.  We want
our intersection graphs to detect when this order is reversed.

\begin{defn} \cite{me3} \label{D:IGstring}
Let $D$ be a chord diagram with $k$ components (oriented line
segments, colored from 1 to $k$) and $n$ chords.  The {\it intersection graph}
$\Gamma(D)$ is the labeled, directed multigraph such that:
\begin{itemize}
    \item $\Gamma(D)$ has a vertex for each chord of $D$.  Each vertex is labeled
by an unordered pair $\{i,j\}$, where $i$ and $j$ are the labels of the components
on which the endpoints of the chord lie.
    \item There is a directed edge from a vertex $v_1$ to a vertex $v_2$ for each
pair $(e_1, e_2)$ where $e_1$ is an endpoint of the chord associated to $v_1$,
$e_2$ is an endpoint of the chord associated to $v_2$, $e_1$ and $e_2$ lie on the
same component of $D$, and the orientation of the component runs from $e_1$ to
$e_2$ (so if the components are all oriented upwards, $e_1$ is below $e_2$).  We
count these edges "mod 2", so if two vertices are connected by two
directed edges with the same direction, the edges cancel each other.  If two vertices
are connected by a directed edge in each direction, we will simply connect them
by an undirected edge.
\end{itemize}
\end{defn}

\begin{figure}
    $$\includegraphics{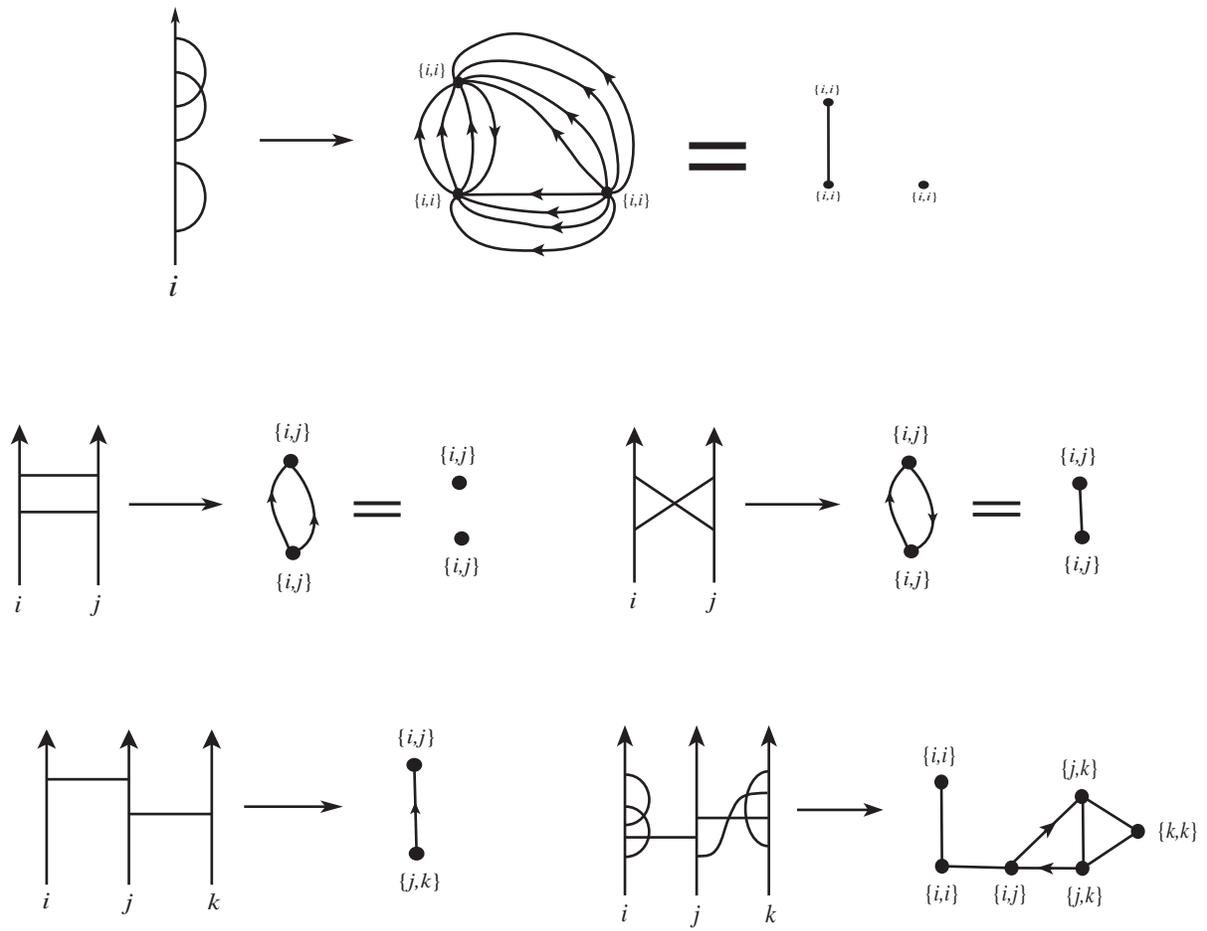}$$
    \caption{Examples of intersection graphs for string links} \label{F:IGstring}
    \end{figure}

Examples of chord diagrams and their associated intersection graphs are given in
Figure~\ref{F:IGstring}.
Note that when the two chords have both endpoints on the same component $i$, our
definition of intersection graph corresponds to the usual intersection graph for
knots.  Our definition also matches our intuition in the case of chord diagrams
of two components, as shown in Figure~\ref{F:IGstring}.

Note also that the total number of directed edges between a vertex $v$ labeled
$\{i,j\}$ and a vertex $w$ labeled $\{l,m\}$ is given by the sum of the number
of occurrences of $i$ in $\{l,m\}$ and the number of occurrences of $j$ in
$\{l,m\}$.  In particular, if a vertex $v$ has a label $\{i,i\}$, this number
will be even (0, 2 or 4).  Since we count directed edges modulo 2, this implies
there is an (uncancelled) directed edge from $v$ to another vertex $w$ if and only if
there is also an (uncancelled) directed edge from $w$ to $v$.  We will say that labeled
directed multigraphs which have this property are {\it semisymmetric}.

\begin{defn} \label{D:semisym}
A directed multigraph G, with each vertex labeled by a pair \{i,j\}, is {\bf
semisymmetric} if for every vertex v labeled \{i,i\}, and any other vertex w, there is a
directed edge from v to w if and only if there is a directed edge from w to v.
\end{defn}

For convenience, we will refer to vertices with labels $\{i,i\}$ (and their associated chords) as {\it
unmarked}, and to other vertices and chords as {\it marked}.  In other words, a marked chord connects
two different components of the chord diagram, whereas an unmarked chord has both endpoints on the same
component.  If we do not need the full structure of the labels, and need only distinguish these two
types of vertices, we will talk about {\it marked}, rather than labeled, trees.  It is understood that
any tree which is labeled as in Definition \ref{D:IGstring} is also a marked tree.

Just as it is useful to consider shares as parts of a chord diagram, it is useful to consider subsets
of the vertices of a tree intersection graph.

\begin{defn} \label{D:bough}
A {\bf bough} of a vertex v in a marked tree T is a connected component of $T\backslash v$
(the graph which results from removing the vertex v and all edges incident to v).  The bough
is called {\bf light} if it contains at most one marked vertex and this vertex (if present)
is adjacent to v in T.  Otherwise the bough is {\bf heavy}.
\end{defn}

In fact, we have shown \cite{me4} that there is a close connection between the boughs of an intersection
graph and shares in the chord diagram.

\begin{lem} \cite{me4} \label{L:share}
Assume T is a marked tree which is the intersection graph for a chord diagram D, and v is a
vertex in T.  Further assume that $T$ has at least one marked vertex.  Then a bough of v is
light if and only if the corresponding chords are a share in D.
\end{lem}

\section{Tree Diagrams with 2 Components} \label{S:2comp}

In this section we will show that the intersection graph determines the (string link) chord diagram
(modulo the 1-term and 4-term relations) in the special case when the intersection graph is a  {\it
trimmed} tree (for diagrams with 2 components).  This generalizes a result of Chmutov, {\it et al.}
\cite{cdl}, who proved that for knots (i.e. string links with one component), tree diagrams are
determined by their intersection graph. 

\begin{defn} \label{D:trimmed}
An intersection graph T for a chord diagram with 2 components is a {\bf trimmed} tree if it is an
(undirected) tree, and there is some vertex v in T, called the {\bf trunk} of T, such that every bough
of v is light (see Definition \ref{D:bough}).
\end{defn}

If $T$ is a trimmed tree with trunk $v$, then every marked vertex in $T$ is adjacent to $v$
(otherwise, $v$ would have a heavy bough).  Moreover, every other vertex in $T$ has at most one heavy
bough; namely, the bough containing $v$.  Also, by Lemma \ref{L:share}, every bough of $v$ corresponds
to a share in $D$ (assuming $T$ has at least one marked vertex).  Our goal in this section is to prove:

\begin{thm} \label{T:2comp}
If $\Gamma(D_1) = \Gamma(D_2) = T$, where $T$ is a trimmed tree, then $D_1$ and $D_2$ are equivalent
modulo the 1-term and 4-term relations.
\end{thm}

We will describe a set of {\it elementary transformations} which will provide an equivalence relation
between chord diagrams with the same intersection graphs, and then show that these transformations can
be achieved using the 1-term and 4-term relations.

\subsection{Elementary Transformations} \label{SS:elementary}

Consider a chord diagram $D$ on two components, whose intersection graph $T$ is a trimmed tree.  Denote
the trunk of $T$ by $v$; we will also use $v$ to denote the corresponding chord in $D$.  Consider the
vertices in $T$ adjacent to $v$ and the corresponding chords in $D$.  These vertices may be both
marked and unmarked:  denote the unmarked vertices $v_1,..., v_n$ and the marked vertices
$y_1,...,y_m$.  Note that if $w$ is any {\it other} vertex of $T$, $w$ will be adjacent to at most one
marked vertex, contained in the same bough of $w$ as $v$.

Since $T$ is a tree, each of $v_i$ and $y_j$ belongs to a different bough of $v$, and so the
corresponding chords belong to non-intersecting shares in $D$.  If $D$ is drawn to minimize crossings
between chords, these chords will cross $v$ in some order; without loss of generality, say the order is
$v_1, v_2,..., v_k, y_1,...,y_m, v_{k+1},..., v_n$.  (Since $T$ is a tree, the $y_j$'s must be grouped
together, because an unmarked chord cannot have both endpoints between two of the $y_j$'s and still
cross $v$.)  Similarly, if $w$ is some other chord of $D$, then the chords crossing $w$ can be put
in order $w_1,...,w_k,y,w_{k+1},..., w_r$, where $y$ corresponds to the heavy bough of $w$ in $T$ (if it
exists), and all the $w_i$'s are unmarked chords corresponding to light boughs in $T$, and hence shares
in $D$.

Now we can define the elementary transformations.

\begin{defn} \label{D:elementary}
The {\bf elementary transformations} of a chord diagram $D$ whose intersection graph is a trimmed tree
are:
\begin{enumerate}
     \item  Permuting the boughs along a chord $w$.  If $w$ is the trunk of $D$, then all the marked
boughs must remain adjacent.  If $w$ is not the trunk, and has a heavy bough, the heavy bough stays
fixed while the other boughs move around it.  Moreover, unmarked boughs cannot be moved from one
component to the other (as this will change the intersection graph).
	    \item  If the trunk of $D$ is marked, reflecting the other marked boughs across the trunk.
\end{enumerate}
Examples of these transformations are shown in Figure~\ref{F:elementary} (A and B are boughs).
    \begin{figure} [h]
    $$\includegraphics{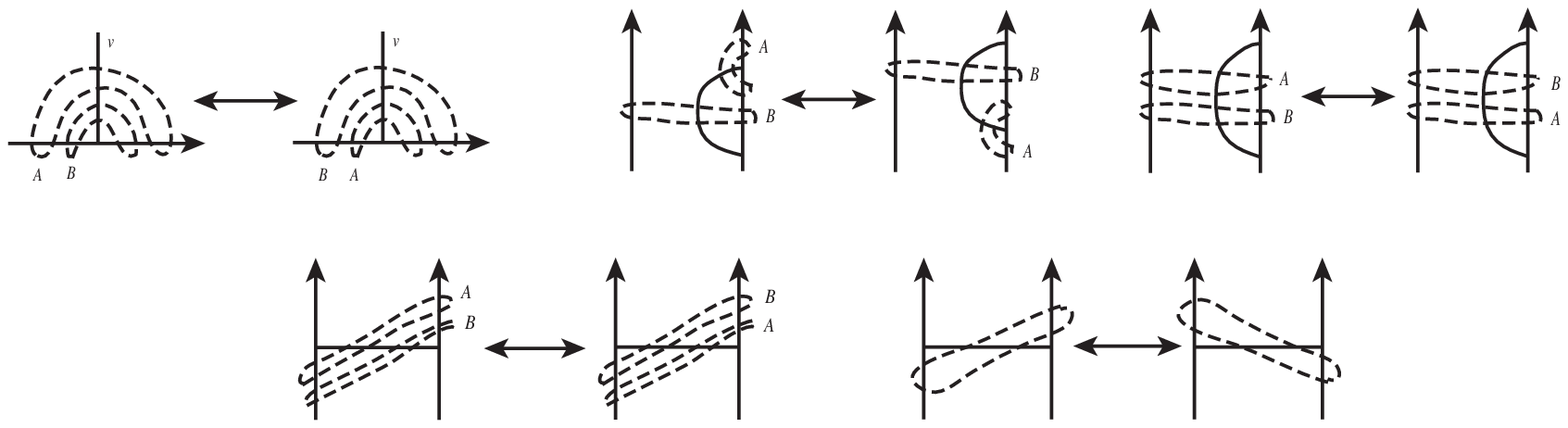}$$
    \caption{Elementary transformations} \label{F:elementary}
    \end{figure}
\end{defn}

\begin{prop} \label{P:elementary}
Let $D_1$ and $D_2$ be chord diagrams such that $\Gamma(D_1)$ and $\Gamma(D_2)$ are trimmed trees.  Then
$\Gamma(D_1) = \Gamma(D_2)$ if and only if $D_1$ can be transformed into $D_2$ via elementary
transformations.
\end{prop}
{\sc Proof:}  It's clear that the elementary transformations have no effect on the intersection graphs,
so we only need to show that if the intersection graphs are the same, then the diagrams are equivalent
modulo the elementary transformations.  Assume that $\Gamma(D_1) = \Gamma(D_2) = T$, and let $v$ be the
trunk of $T$.  The order of the boughs along $v$ in $D_1$ differs from the order along $D_2$ by a
permutation which keeps all the marked boughs adjacent (since they must be adjacent in both diagrams),
and which keeps unmarked chords on the same component; this can be achieved by an elementary
tranformation of type 1.  If $v$ is marked, it is possible that its marked boughs in $D_1$ differ from
the marked boughs in $D_2$ by the direction of their slant - this can be corrected by an elementary
transformation of type 2.  We can now consider the vertices adjacent to $v$.  Once again, we can
rearrange the order of the boughs by elementary transformations of type 1; moreover, we can consider
the bough containing $v$ (the only possible heavy bough) to be fixed during this process, so we will
not effect the results of our previous moves.  Since these chords can have at most one marked bough
(the one containing $v$), we will never need to use transformations of type 2.  Continuing inductively,
moving on to vertices farther and farther from the trunk, we can rearrange the chords of $D_1$ until
$D_1 = D_2$. $\Box$

\subsection{Intersection Graph Conjecture} \label{SS:IGC}

In this section we will prove Theorem~\ref{T:2comp}.  Our proof is modeled on the proof for knots given
by Chmutov {\it et al} \cite{cdl}.  Using Proposition~\ref{P:elementary}, it is sufficient to show
the following:

\begin{prop} \label{P:equivalence}
If $D_1$ and $D_2$ are chord diagrams on two components which differ by an elementary transformation,
then $D_1$ and $D_2$ are equivalent modulo the 1-term and 4-term relations.
\end{prop}

Before we begin to prove this proposition, we will state a few useful facts.

\begin{lem} \label{L:gen4term}
(Generalized 4-term relation \cite{cdl})  For any share and chord, the following relation holds (modulo
the usual 4-term relation):
$$\includegraphics{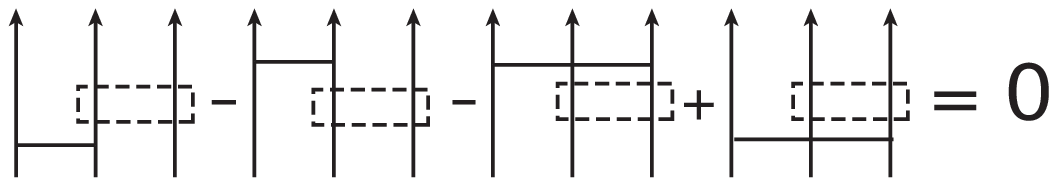}$$
\end{lem}
{\sc Proof:}  Write down the 4-term relations for the given chord and every chord in the given share,
and add together all these relations.  Except for the four terms of the generalized 4-term relation,
every term will appear twice with opposite sign, and so cancel.  We are left with the generalized 4-term
relation. $\Box$

\begin{cor} \label{C:simple}
\cite{cdl,do} Modulo the 1-term and 4-term relations, we have the following relations among chord
diagrams:
$$\includegraphics{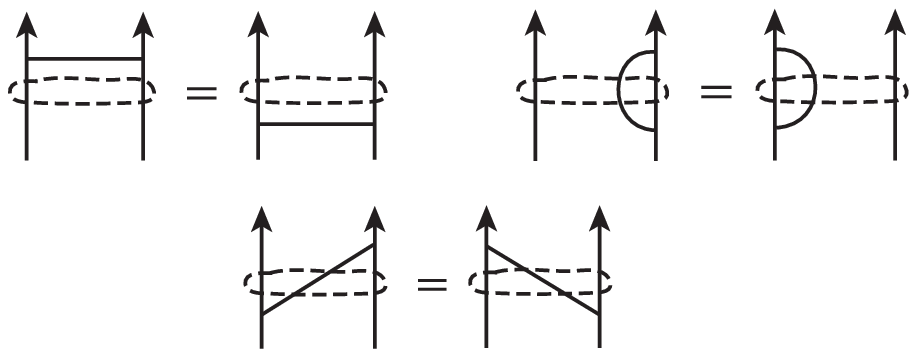}$$
\end{cor}
{\sc Proof:}  For each pair of diagrams, keeping the shares fixed, consider all diagrams formed from
the shares and a single chord.  Writing down all the generalized 4-term relations among these diagrams,
and cancelling terms using the 1-term relation, we quickly see that the three relations are
equivalent.  Dogolazky \cite{do} proved the first relation holds. $\Box$

{\sc Remark:}  Chmutov {\it et al.} \cite{cdl} proved the result by solving the system of equations
arising from the 4-term relations.  Dogolazky \cite{do} pointed out that this proof only works if the
space $B_n^k$ has no elements of order 2, and that this is not true in general.  However, he provided a
more complex proof to show the result is true in general.

Our proof of Proposition~\ref{P:equivalence} will consider several cases.  Since any permutation can be
decomposed into transpositions, it will suffice to consider a few transpositions of boughs.

\begin{lem} \label{L:permute}
The four pairs of diagrams below are equivalent modulo the 1-term and 4-term relations.  The diagram
consisting of the shares $A$ and $B$ and the chord $v$ is a tree diagram; there are no restrictions on
the share $S$.  There are no other chords with endpoints on the indicated arcs, although there may be
other chords which intersect the chords shown.

(1) $$\includegraphics{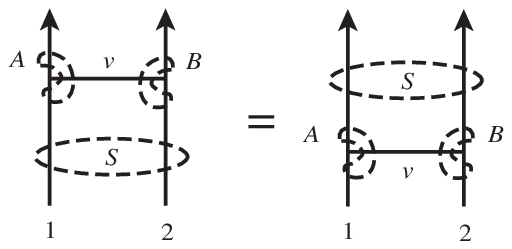}$$
(2) $$\includegraphics{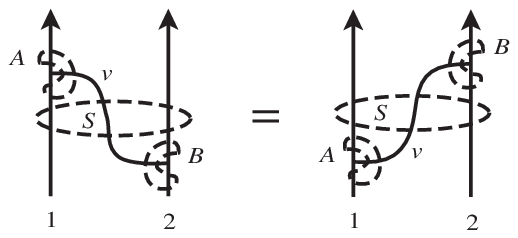}$$
(3) $$\includegraphics{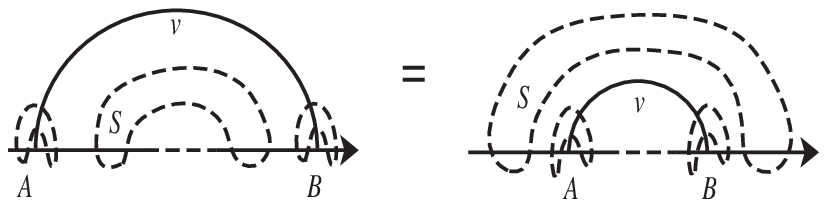}$$
(4) $$\includegraphics{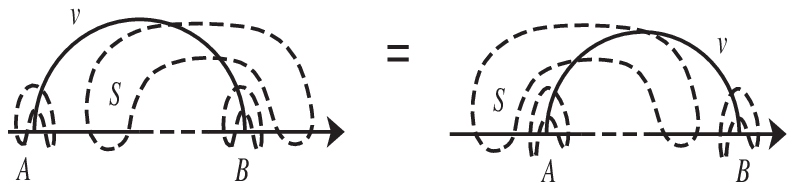}$$
\end{lem}

Our proof will proceed by induction on the total number of chords in the shares $A$ and $B$; we
will call this the {\it complexity} of the diagram $D$, denoted $c(D)$.  Our proof of
Lemma~\ref{L:permute} will require two additional lemmas.  The first lemma is essentially identical to
a lemma in \cite{cdl}; we provide a proof here for clarity and to show the (very slight) modifications
needed for diagrams on string links.

\begin{lem} \label{L:rewrite}
Suppose that Lemma \ref{L:permute} holds for any diagram with complexity less than m.  Then:
$$\includegraphics{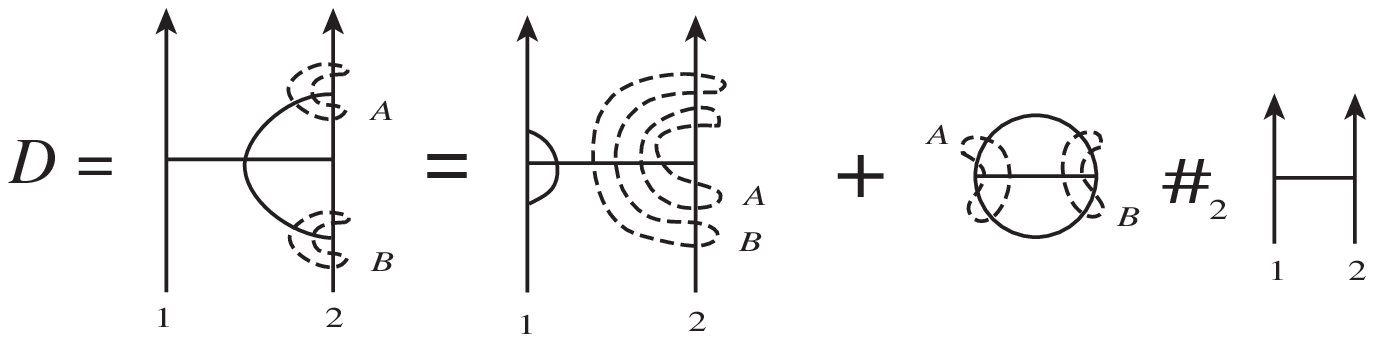}$$
where the total number of chords in the shares A and B is $\leq m$.
\end{lem}
{\sc Proof of Lemma \ref{L:rewrite}:}  By the generalized 4-term relation, we have:
$$\includegraphics{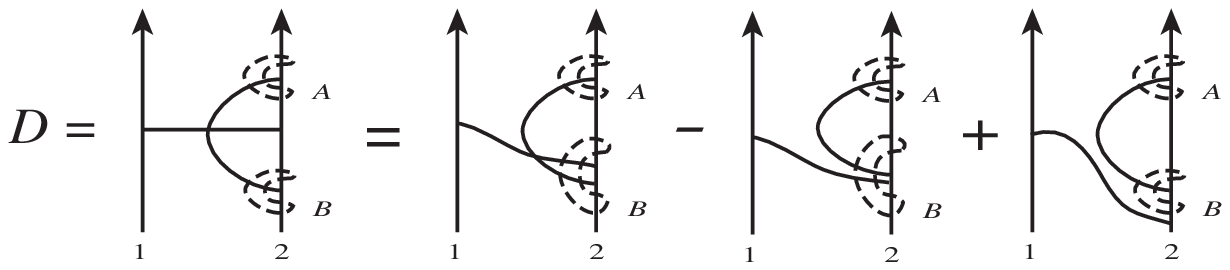}$$
Using the generalized 4-term relation again, and rewriting the last term as a product, this equals:
$$\includegraphics{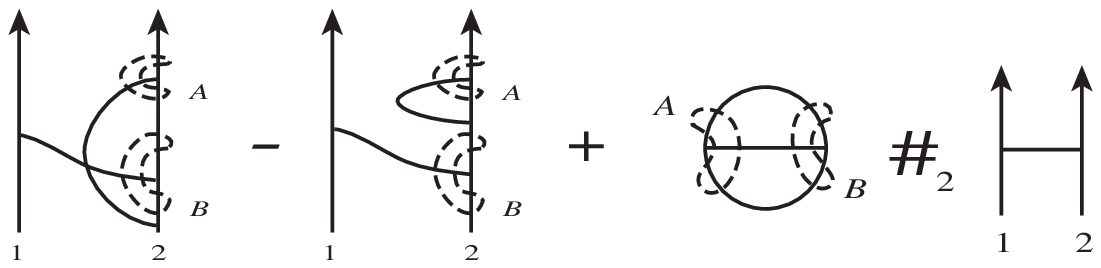}$$
Since the number of chords in the share $A$ is less than $m$ (unless $B$ is empty, which is a trivial
case) we can apply Lemma \ref{L:permute} and permute the boughs of the first diagram.  And since the
multiplication $\#_2$ is well-defined, we can rewrite the second diagram to get:
$$\includegraphics{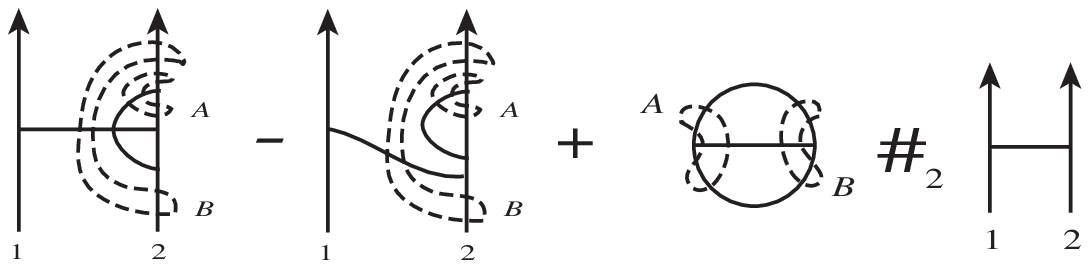}$$
Applying the usual 4-term relation gives:
$$\includegraphics{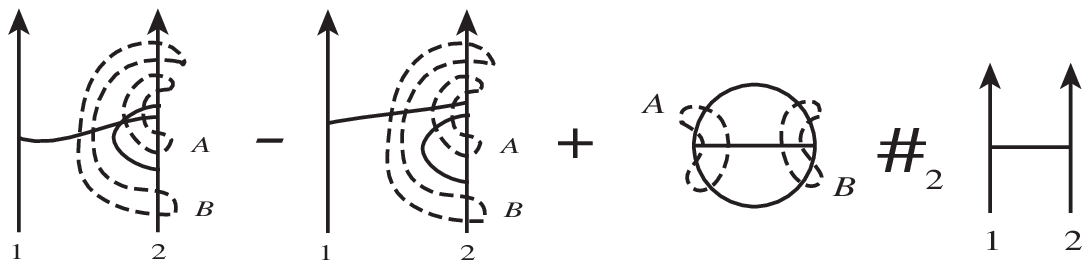}$$
Using the generalized 4-term relation once more we get:
$$\includegraphics{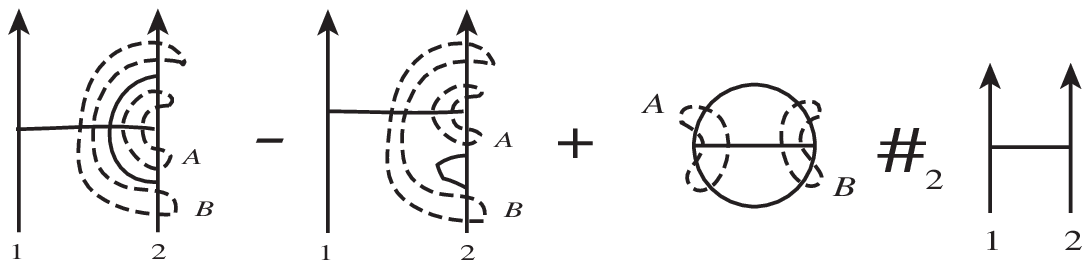}$$
The second diagram vanishes by the 1-term relation.  We can apply Corollary \ref{C:simple} to the first
diagram to obtain:
$$\includegraphics{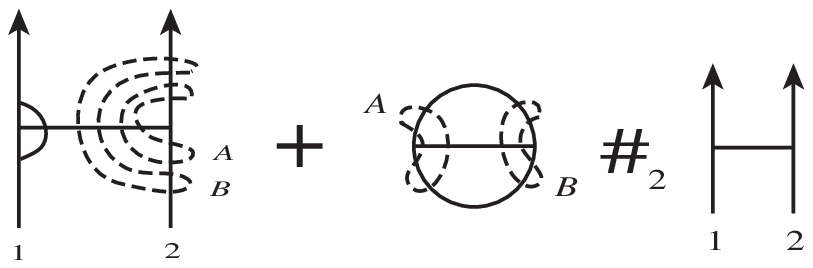}$$
which completes the proof.  $\Box$

\begin{cor} \label{C:deconstruct}
Assume that Lemma \ref{L:permute} holds for diagrams of complexity less than m.  Then:
$$\includegraphics{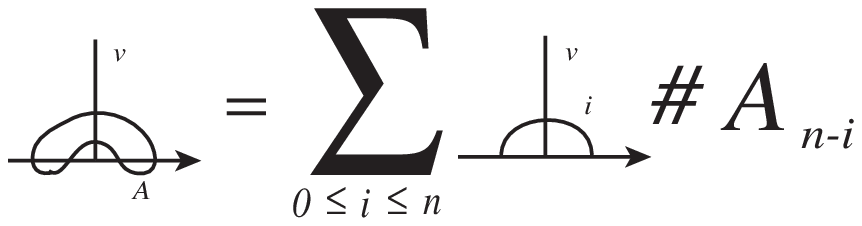}$$
where the share A is a tree diagram, $n \leq m$ is the number of chords in the share A, the chord
labeled i represents i disjoint chords along the chord v, and $A_{n-i}$ is the sum of all diagrams of
degree n-i obtained by taking the connect sum of disjoint boughs of A.  $A_0$ is the empty diagram,
with no chords.
\end{cor}
{\sc Proof:}  We will prove this corollary by induction on $n$.  In the base case, $n = 1$, and $A$ is
just a single chord.  Then $A_{n-i}$ is either $A_0$ (the empty diagram) or $A_1 = A$.  But $A_1 = 0$,
by the one-term relation.  So there is only one term on the right-hand side of the equation, when $i =
1$, which gives us the original diagram.

Now assume that the corollary holds when there are fewer than $n$ chords in the share, and assume $A$
has $n$ chords.  We have two cases:  either $A$ is connected (i.e. $A$ is a single bough of the chord
$v$), or $A$ is made up of several boughs of $v$.  We will first consider the case when $A$ is the
union of several boughs.  Then $A$ can be divided into two shares $B$ and $C$, each of which contains
at least one bough of $v$.  Say that $B$ has $m$ chords and $C$ has $k$ chords, so $0 < m,k < n$. 
Then, by our inductive hypothesis, we can apply the corollary to the shares $B$ and $C$, as follows:
$$\includegraphics{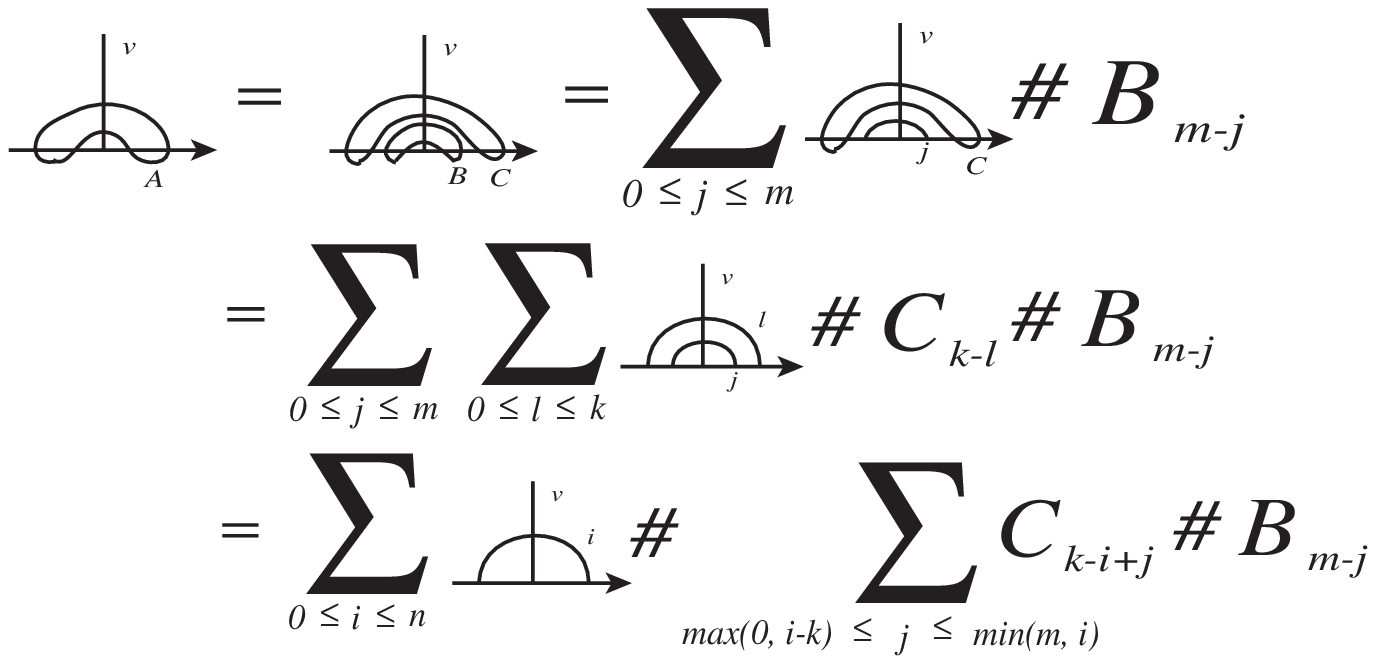}$$
Since each term of $A_{n-i}$ is simply the product (via connect sum) for some $j \leq i$ of a term of
$B_{m-j}$ and a term of $C_{k-i+j}$, we are left with the desired sum.

In the case when $A$ is a single bough of $v$, we can apply Lemma \ref{L:rewrite} (since $n \leq m$,
so the total number of chords in the shares $B$ and $C$ below is less than $m$):
$$\includegraphics{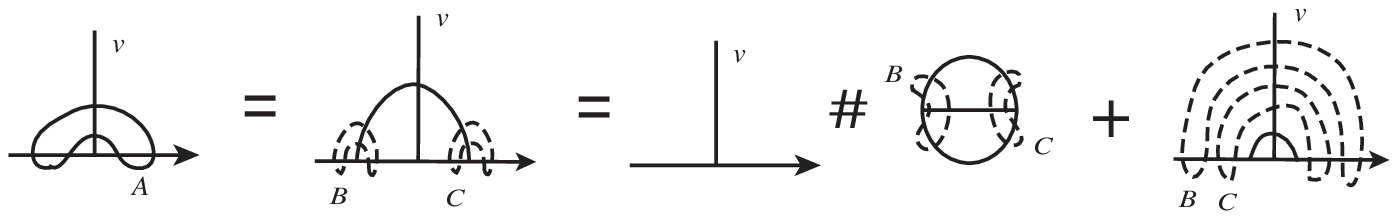}$$
The first term is the term for $i=0$; the other term can be decomposed into the terms for $i > 0$ as in
the previous case.  This completes the proof of the corollary. $\Box$

\begin{lem} \label{L:center}
The diagram $J_n$ shown below is in the center of $B_2$ for $n \geq 0$.
$$\includegraphics{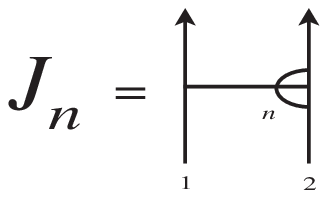}$$
\end{lem}
{\sc Proof:}  Dogolazky \cite{do} showed that certain linear combinations of chord diagrams $I_n$
(represented as a {\it unitrivalent diagram}) are in the center of $B_2$ for $n \geq 1$.  A discussion
of the isomorphism between the spaces of chord diagrams and unitrivalent diagrams can be found in
\cite{bn}; we will not reproduce it here, since it is not required for the remainder of the paper.  It
is an easy exercise in induction to show that $J_n = (-1)^n I_{n+1}$, so $J_n$ is also in the center.
$\Box$
\\
\\
\noindent{\sc Proof of Lemma \ref{L:permute}:}  We will start by proving relation (1); the proofs of
the other parts are similar.  We are inducting on the complexity of the chord diagram (i.e. the number
of chords in the shares $A$ and $B$).  The base case is when the complexity is 0; in this case the
result is given by Corollary \ref{C:simple}.

For our inductive step, we assume the diagram has complexity $c(D) = m$, and assume the lemma holds for
all diagrams with complexity less than $m$.  Denote the number of chords in the shares $A$ and $B$ by
$a$ and $b$ respectively.  Since $a$ and $b$ are both at most $m$, we can apply Corollary
\ref{C:deconstruct} to the diagram on the left-hand side to get:
$$\includegraphics{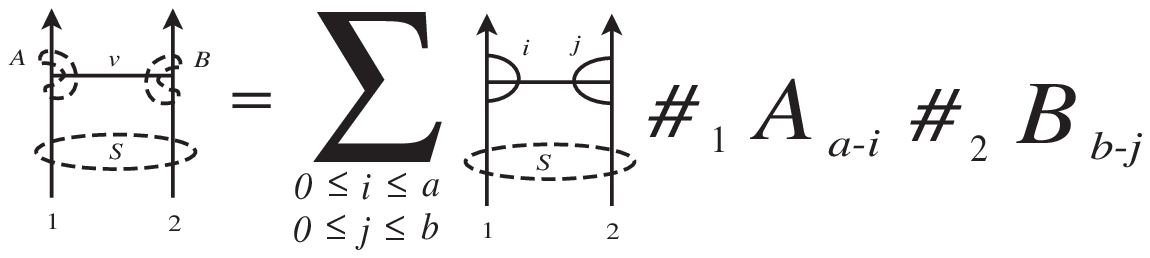}$$
Here, a chord labeled, for example, $i$ denotes a share consisting of $i$ disjoint chords.  So it
suffices to consider the case when the shares $A$ and $B$ each consist of collections of disjoint
chords.  But by Corollary \ref{C:simple}, this means the diagram consisting of the shares $A$ and $B$
and the chord between them is simply $J_{a+b}$, which is in the center of $B_2$ by Lemma \ref{L:center}.

To prove the other three relations, notice that relation (4) is the same as relation (1), except that
the oriented line segments belong to the same component.  Relation (2) is also the same as (1),
except that we have reversed the orientation of one line segment.  But, by the antisymmetry relation,
this only changes the diagrams by a sign, so the equality still holds.  Finally, relation (3) is the
same as (2), when the line segments belong to the same component.  This completes the proof. $\Box$

So far, we have only shown that we can permute neighbouring marked or unmarked boughs.  The next case
we need to consider is moving an umarked bough from one end of an unmarked chord to the other, possibly
across marked chords.

\begin{lem} \label{L:endtoend}
The two diagrams below are equivalent modulo the 1-term and 4-term relations.  No other chords have
endpoints on the solid line segments, though there may be chords (including chords to other components)
with endpoints on the dashed segments.
$$\includegraphics{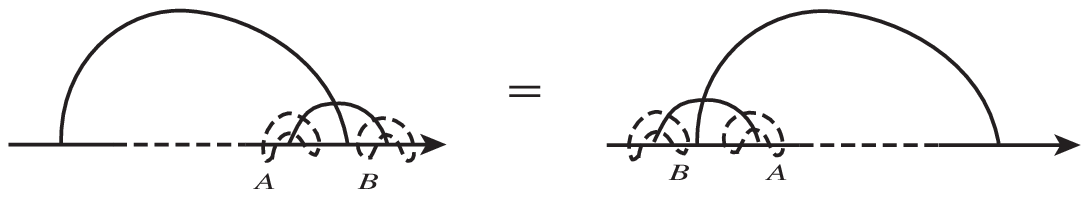}$$
\end{lem}
{\sc Proof:}  As before, our proof is by induction on the complexity of the diagram - the number of
chords in the shares $A$ and $B$.  If the complexity is 0, then $A$ and $B$ are empty, and the relation
follows from Corollary \ref{C:simple}.  So we may assume the relation holds for diagrams with
complexity less than $m$. 

Assume the diagrams in the relation have complexity $m$.  If we apply Lemma \ref{L:rewrite} to the
diagram on the left-hand side, we obtain:
$$\includegraphics{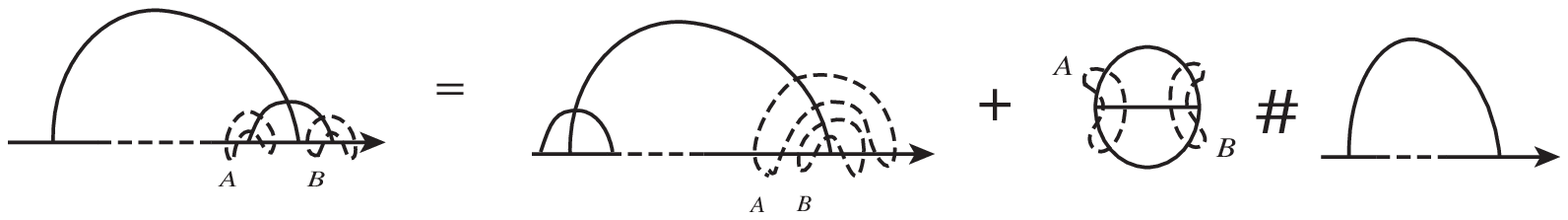}$$
Applying Lemma \ref{L:rewrite} to the diagram on the right-hand side gives almost the same result,
except that the shares $A$ and $B$ (and the single chord) are at the other ends of the chord.  But
since these shares have fewer than $m$ chords, we may apply our inductive hypothesis to move them to
the other end as needed.  This completes the proof. $\Box$

All that remains is to show that diagrams which are equivalent modulo the second elementary
transformation (flipping the marked chords across the trunk of the tree) are also equivalent modulo the
1-term and 4-term relations.  But, in fact, we have already done this - this is part (2) of Lemma
\ref{L:permute}.  This completes the proof of Proposition \ref{P:equivalence}.  Together, Propositions
\ref{P:elementary} and \ref{P:equivalence} prove Theorem \ref{T:2comp}.

\begin{figure} [h]
    $$\includegraphics{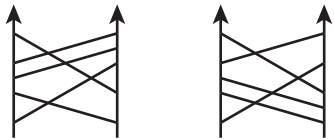}$$
    \caption{Unequal tree diagrams with the same intersection graph.} \label{F:untrimmed}
    \end{figure}

{\sc Remark:}  Perhaps surprisingly, if the intersection graph is an untrimmed tree it does {\it not}
determine the diagram modulo the 1-term and 4-term relations (for diagrams on two components).  A
counterexample was given by Dogolazky \cite{do}, who showed by a computer calculation that the two
diagrams in Figure~\ref{F:untrimmed} are not equal in $B_2^5$.

\section{Tree Diagrams with more than 2 Components}
\label{S:>2comp}

The case of tree diagrams with more than two components is actually simpler than the case of two
components, since there are greater restrictions on the trees which can be realized as intersection
graphs.  The following theorem was proven in \cite{me4}:

\begin{thm} \label{T:treeclass}
Let T be a connected, labeled directed tree with n colors (so each vertex has a label \{i, j\}, where $1
\leq i, j\leq n$).  Let $m_{i,j}$ denote the number of vertices with label \{i, j\}.  Then T is an
intersection graph for a connected chord diagram on n components if and only if the following
conditions are met (possibly after relabeling the tree by a permutation of 1,..., n):
\begin{enumerate}
   \item The labels of adjacent vertices must have at least one color in common.
   \item T is semisymmetric (see Definition \ref{D:semisym}).
   \item If v has label \{i,j\} and w has label \{i,k\}, where i, j and k are all distinct, then there is
a directed edge between v and w.
   \item $m_{i,j} = 0$ if $|i-j| > 1$.
	  \item $m_{i, i+1} = 1$ for $2\leq i\leq n-2$, $m_{1,2}\geq 1$ and $m_{n-1,n} \geq 1$.
   \item No two marked vertices are connected by a path of undirected edges.
\end{enumerate}
\end{thm}

Let $T$ be a tree satisfying the six conditions above, so it is the intersection graph for some
chord diagram $D$.  Our goal is to reconstruct $D$ from $T$, at least modulo the 1-term and 4-term
relations.  Consider the graph $F$ which results from removing all the directed edges in
$T$.  By the last condition above (and the fact that $T$ was connected), each
component of this diagram will contain a single marked vertex.  Exactly one of these marked vertices
will be labeled $\{i,i+1\}$ unless $i=1$ or $i=n-1$, in which cases there may be more than one.  So each
component of $F$ is the intersection graph for a trimmed tree diagram on two components, and by Theorem
\ref{T:2comp} this diagram is determined by its intersection graph (modulo the 1-term and 4-term
relations).  All that remains is to arrange these diagrams to form $D$.  This arrangement is determined
by the directions of the directed edges of $T$, except for the components with marked chords labeled
$\{1,2\}$ or $\{n-1,n\}$.  For example, if $n > 4$ and there is a directed edge from the (unique)
vertex of $T$ labeled $\{2,3\}$ to the (unique) vertex labeled $\{3,4\}$, then the chords corresponding
to the component of $F$ containing the first marked chord will be placed below the chords corresponding
to the component of $F$ containing the second marked chord.  The chords at the ends need to be treated
a little more carefully.  The components with a chord labeled $\{1,2\}$ naturally fall into two groups:
those which lie below the (unique) chord labeled $\{2,3\}$, and those which lie above this chord. 
Within these two groups, however, the order of the components can be rearranged without affecting the
intersection graph.  However, these rearrangements are all permutations of the components, which by
Lemma \ref{L:permute} do not change the diagram modulo the 1-term and 4-term relations.  The same is
true for the components with a chord labeled $\{n-1,n\}$.  So, modulo the 1-term and 4-term relations,
we can reconstruct the diagram $D$ from the intersection graph $T$.  We can conclude:

\begin{thm} \label{T:>2comp}
If $D_1$ and $D_2$ are connected tree diagrams on n components ($n > 2$) which have the same
intersection graph, then $D_1$ and $D_2$ are equivalent in $B_n$.
\end{thm}

Note that, unlike for diagrams on 2 components, we do not need to restrict to a smaller class of tree
diagrams, because of the greater constraints imposed by Theorem \ref{T:treeclass}.

\section{Questions} \label{S:questions}

We noted at the end of Section \ref{S:2comp} that, in general, tree diagrams on two components are not
determined by their intersection graphs.  However, the example found by Dogolazky \cite{do} lies in the
torsion subgroup of $B_2^5$ - the difference $\Delta w$ of the two diagrams in Figure \ref{F:untrimmed}
has order 2.

\begin{quest} \label{Q:quest1}
Do all counterexamples to the Intersection Graph Conjecture for tree diagrams have finite order?  In
other words, is the conjecture true if we look at the quotient of $B_2$ by its torsion subgroup?
\end{quest}

Stanford has shown that, while $\Delta w$ is nontrivial in $B_2^5$, any realization of it as a
difference of two singular string links {\it is} trivial modulo the "topological" 1-term and 4-term
relations.  So one could also ask:

\begin{quest} \label{Q:quest2}
Is any counterexample to the Intersection Graph Conjecture for tree diagrams trivial when realized as
a linear combination of singular string links?
\end{quest}

Given that intersection graphs characterize chord diagrams for trimmed trees, it would be interesting
to use this to find a basis for the space of trimmed tree diagrams.

\begin{quest} \label{Q:quest4}
What is the dimension of the space of chord diagrams spanned by diagrams whose intersection graphs are
trimmed trees?  What is a basis for this space?
\end{quest}

And, of course, it is natural to ask whether we can move beyond trees.

\begin{quest} \label{Q:quest3}
To what extent do intersection graphs determine chord diagrams in general?
\end{quest}

\small

\normalsize

\end{document}